\documentstyle[righttag,12pt]{amsart}

\newtheorem{thm}{Theorem}[section]
\newtheorem{prop}[thm]{Proposition}
\newtheorem{cor}[thm]{Corollary}
\newtheorem{lem}[thm]{Lemma}

\newtheorem{remark}[thm]{Remark}

\newcommand{\cn}{{\Bbb C}^n}
\newcommand{\p}{\partial}

\newcommand{\ws}{W^s(\Omega)}
\newcommand{\wsq}{W^s_{(0,q)}(\Omega)}

\def\pd#1#2{\frac{\p #1}{\p #2}}
\def\e#1#2{\varepsilon^{#1}_{#2}}
\def\vt{
\vartheta}

\def\dom{\text{dom\,}}
\def\rg{\text{range\,}}
\def\ker{\text{ker\,}}
\def\ss{\subseteq}
\def\dbar{\overline{\partial}}
\def\dbars{(\dbar,s)}

\def\la{\langle}
\def\ra{\rangle}
\def\K{{\cal K}}

\def\CP{\gamma_\alpha}
\def\CPp{\gamma_{\alpha'}}
\def\CPb{\gamma_\beta}
\def\bvp{boundary value problem}
\def\po{\p_{x_0}}
\def\bn#1#2{{\textstyle \binom{#1}{#2}}}

\begin{document}
\title[The $\dbar$-Neumann problem in the 
Sobolev
topology]{The $\dbar$-Neumann problem in the Sobolev  
topology}

\author[L. 
Fontana]{Luigi Fontana}
\address{\hskip-\parindent
Luigi Fontana, Dipartimento di Matematica, Via Saldini 50,
Universit\`a di Milano\\ 
20133 Milano (Italy)}
\email{fontana@@vmimat.mat.unimi.it}
\author[S. G. Krantz]{Steven G. Krantz}\thanks{Krantz's research was supported 
in part by Grant DMS-9531967
from the National Science Foundation.  Research
at MSRI is supported by NSF Grant DMS-9022140.}
\address{\hskip-\parindent Steven G. Krantz,
Mathematical Sciences Research Institute {\normalfont and}
Department of Mathematics, Washington University, St.
Louis, MO 63130 (U.S.A.)}
\email{sk@@artsci.wustl.edu}
\author[M. M. Peloso]{Marco M. Peloso}
\address{\hskip-\parindent
Marco M. Peloso, Dipartimento di Matematica, Politecnico di Torino, 10129 
Torino (Italy)}
\email{peloso@@polito.it}


\begin{abstract}
We study the $\overline{\partial}$-Neumann problem
using the Sobolev space inner product.  We show that the problem
can be solved on any smoothly bounded, pseudoconvex domain.
We further formulate estimates and the basic results of
a Sobolev Hodge theory.
\end{abstract}

\maketitle 

\section{Introduction}
Let $\Omega$ be a smoothly bounded domain in $\cn$. 
We write the coordinates $z_j =x_j +ix_{j+n}$, $j=1,\dots,n$, 
and the standard basis of vector fields $D_k := \p/\p x_k$, for
$k=1,\dots,2n$. 
For $s$ a
non-negative integer we define the Sobolev 
inner product $\la \cdot,\cdot\ra_s$ to be
\begin{equation}\label{Sobolev}
\la f,g\ra_s := \sum_{|\alpha|\le s} \CP \int_\Omega D^\alpha f 
\overline{D^\alpha g} .
\end{equation}
Here, and throughout the paper, we use $D^\alpha$ to denote
the $\alpha$-order derivative, where $\alpha$ is a multi-index
and we are using standard multi-index notation.  Moreover, 
$\CP:=|\alpha|!/\alpha!$ denotes the polynomial coefficient. 
[The naturality of this choice of the Sobolev inner 
product will be
pointed out and discussed below.]

We define the Sobolev
space $W^s (\Omega)$ to be the closure of
$C^\infty (\bar \Omega)
$ with respect to the above inner product. 
We denote by $\wsq$ the space of $(0,q)$ forms 
whose coefficients are
in $\ws$. If $\phi=\sum_{|J|=q}\phi_J d\bar z^J$ and
$\psi=\sum_{|J|=q}\psi_Jd\bar z^J$, then the inner product in $\wsq$ 
is defined by
$$
\la \phi,\psi\ra_s := \sum_{|J|=q}\sum_
{|\alpha|\le s} \CP
\int_\Omega D^\alpha \phi_J \overline{D^\alpha \psi_J} ,
$$
where we use the standard notation $J$ to denote a $q$-vector with
increasing entries, and $\alpha$ to denote a multi-index.
[Note that the inner product of forms of different degrees
is defined to be 0.]

For a $(0,q)$ form $\phi=\sum_{|J|=q}\phi_J d\bar z^J$ 
with $C^\infty$ coefficients, the operator $\dbar$
is defined by
\begin{equation}\label{d-bar}
\dbar\phi :=\sum_{|K|=q+1} \sum_{kJ}\e{K}{kJ} \pd{\phi_J}{\bar z_k}  
d\bar z^K ,
\end{equation}
where $\e{K}{kJ}$  equals the sign of the permutation $kJ\mapsto K$
if $\{ k\}\cup J=K$ as sets, and is $0$ otherwise.  We continue
to use $\overline{\partial}$ to denote 
its closure in the $W^s$ topology.  In this way,
for each integer $q=0,1,\dots,n$, we obtain an 
unbounded, densely
defined, closed operator
$$
\dbar:\wsq\rightarrow W^s_{(0,q+1)}(\Omega) .
$$
Thus, in particular, $\ker \dbar$ is a closed subspace in $\wsq$.
Sometimes we shall use the notation $\dbar_{(0,q)}$ to stress the
fact that the
operator
 $\dbar$ is acting on $(0,q)$ forms.

Consider now the 
$\ws$-Hilbert space adjoint $\dbar^*$ of $\dbar$.
We want to study the
boundary value problem
\begin{equation}\label{dbar-neumann}
\begin{cases}
(\dbar\dbar^* +\dbar^* \dbar)u=f &\text{ on }\Omega\\
u,\, \dbar u\in \dom \dbar^* \, , &   
\end{cases} 
\end{equation}
where $f$ is a given $(0,q)$ form.
When appropriate, we shall refer to this problem as {\bf (3,s)}
in order to emphasize that the topology is coming from
the $W^s$ inner product.
The condition that $u$ and $\overline{\partial} u$ lie in the
domain of $\overline{\partial}^*$ leads to the
{\em $\overline{\partial}$-Neumann 
$s$-order boundary conditions}. 
We shall refer below to the
$(\dbar,s)$-Neumann conditions, and the $\dbars$-Neumann problem.
Notice that if the Hilbert space under consideration is $L^2
(\Omega)$ (that is, $s=0$) with respect to the Lebesgue measure, 
then the problem
{\bf (3,s)} reduces to the classical $\dbar$-Neumann problem. 

J.\ J.\ Kohn solved the $\dbar
$-Neumann ($= (\dbar,0)$-Neumann)
problem in a series of papers in 1963-4 (see 
\cite{Folland-Kohn} and references therein). This work has proved 
important in the theory of partial differential equations, in
geometry, and in function theory. Recent  work of Christ
\cite{Christ} has shown that the 
{\em canonical solution\/}---the solution 
that is minimal in $L^2$ norm---that arises from
Kohn's work 
in the $L^2$ topology is not as well behaved as one 
might have hoped.
The program presented in this paper
endeavors to seek other canonical solutions that
may serve when Kohn's solution will not. This work is also
interesting from the point of view of partial differential
equations---particularly boundary value problems---and in the study 
of the energy integral in geometry.  
We mention that 
H.\  Boas \cite{Boas1} and \cite{Boas2} studied properties
and regularity of the Hilbert space orthogonal projection 
of $\ws$ onto the subspace consisting of the holomorphic functions.

The present paper is the first of a series of papers that we devote
to the study of the $\dbars$-Neumann problem. We begin by showing
that problem {\bf (3,s)} 
can always be solved on any smoothly bounded
pseudoconvex domain $\Omega$.  
This result does  not depend on the particular choice of Sobolev 
inner product.
Then we investigate the
$\dbars$-Neumann problem more closely by 
determining a description of
the Hilbert space adjoint $\dbar^*$
 of $\dbar$, and the boundary
conditions arising from requiring that $u$ and $\dbar u$ belong to
$\dom\dbar^*$. 
While doing this we use the particular choice of the inner product
(\ref{Sobolev}) to obtain reasonably clean equations and formulas.
We then conclude with some remarks about what lies ahead. In a
forthcoming 
paper we give estimates for the above problem in the special case of 
a strongly pseudoconvex domain, and with $s=1$.  The foundations for
the present work, 
studied in the real variable context of the de Rham
complex, were laid in the papers \cite{FKP1}, \cite{FKP2}.

We thank H. Boas for making several useful 
remarks and comments on an
earlier version of this paper.  We also thank the referee for making
helpful suggestions.  Work of the second author at MSRI
was supported by NSF Grant DMS-9022140.

\section{Solvability of the $\dbars$-Neumann problem}
The aim of the present section is to prove the following theorem.

\begin{thm}\label{solvability} \sl
Let $\Omega$ be a smoothly bounded pseudoconvex domain in $\cn$.
Let $s,q$ be positive integers, $0<q\le n$. 
Let $f\in\wsq$. Then there exists a unique
$u\in\wsq$ that solves the $\dbars$-Neumann problem
$$
\begin{cases}
(\dbar\dbar^* +\dbar^* \dbar)u=f & \text{ on }\Omega\\
u,\, \dbar u\in\dom \dbar^* \, . &
\end{cases} 
$$
Moreover, there exists a
 constant $c>0$, 
independent of $f$, such that
$$
\|u\|_s \le c\|f\|_s.
$$
\end{thm}
\begin{pf}

The proof is in two steps.  
In the first, we  rely heavily on Kohn's estimates \cite{Kohn}, to
construct and estimate the
{\em canonical solutions\/} in $W^s$
to the equations $\dbar u=f$ and
$\dbar^* v=g$.  In the second step we prove the solvabilty of the
$(\dbar,s)$-Neumann problem.  In the course of the proof, by {\em
orthogonal\/} we shall always mean orthogonality in the $W^s$ inner 
product. 

By (3.21) in \cite{Kohn}, since the $\dbar$-cohomology is trivial on
a pseudoconvex domain $\Omega\ss {\Bbb C}^n$, we have that
$\rg \dbar_{(0,q-1)} = \ker \dbar_{(0,q)}$. 
This equality implies that $\rg\dbar_{(0,q-1)}$ is closed in $\wsq$.
Now Lemma 4.1.1 in \cite{HO}, applied with $F=\rg\dbar_{(0,q-1)}$,
gives that 
$$
\| f\|_s \le c \|\dbar^*_{(0,q)} f\|_s 
$$
for all $f\in \rg\dbar_{(0,q-1)}\cap\dom\dbar^*_{(0,q)}$.  This in
turn, by Lemma 4.1.2 in \cite{HO}, implies that for all $v$ in the
orthogonal complement of $\ker\dbar_{(0,q-1)}$, i.e.  in the closure 
of $\rg\dbar^*_{(0,q)}$, there exists $f\in\dom\dbar^*_{(0,q)}$ such
that 
$\dbar^*_{(0,q)} f=v$. 
Hence,
$\rg\dbar_{(0,q)}^*$ is closed as well, and therefore 
we have the estimate $\| f\|_s \le C\|\dbar f\|_s$ for all
$f\in\rg\dbar^*_{(0,q)}\cap\dom\dbar_{(0,q-1)}$. 
Moreover, we have the strong orthogonal decomposition
$$
\wsq=\rg\dbar_{(0,q+1)}^* \oplus\rg\dbar_{(
0,q-1)} .
$$
Now, given any $g\in\wsq$, with $\dbar_{(0,q)} g=0$, i.e.
$g\in\rg\dbar_{(0,q-1)}$, we can find $v\in\dom\dbar_{(0,q-1)}$, 
orthogonal to $\ker\dbar_{(0,q-1)}$, such that $\dbar v=g$, and
we have the estimate
$$
\| v\|_s \le c_s \|g\|_s .
$$
We can apply the same argument to the $\dbar^*$-equation, i.e. given 
any $f$ with $\dbar^*_{(0,q)}f=0$, we can find $u$ orthogonal to
$\ker\dbar^*_{(0,q+1)}$ such that $\dbar^*_{(0,q+1)} u=f$, with the
estimate 
$$
\| u\|_s \le c_s \|f\|_s .
$$
We shall call such solutions $u$ and $v$ the $s$-{\em canonical}
solution to the $\dbar$ and $\dbar^*$ equation, respectively.
  
We now establish the solvability of the $\dbars$-Neumann
problem.
We shall suppress the subscripts on the operators
$\dbar$ and $\dbar^*$ (used to denote the  space of forms
that  is being acted upon),  since this will be clear from context.  
Let $f\in\wsq$. Then $f$ can be uniquely
written as $f=f_1+f_2$ with $f_1 \in \rg\dbar$ and
$f_2\in\rg \dbar^*$. Let $g_1,g_2$ be
 the canonical solution of
$\dbar g_1=f_1$, and
$\dbar^* g_2=f_2$, respectively.
Since $g_1\perp \ker \dbar$ we have that $g_1\in\rg\dbar^*$, and
therefore $\dbar^* g_1=0$. Analogously, 
$g_2\in\rg\dbar$ and $\dbar g_2=0$. 
 
Thus we can canonically select
$u_1,u_2$ such that $\dbar^* u_1=g_1$ and $\dbar u_2 =g_2$.
Setting $u=u_1+u_2$ we obtain that
$$
(\dbar\dbar^* +\dbar^* \dbar)u=f,
$$
and the desired estimate follows from the corresponding ones for 
$\dbar$ and $\dbar^*$:
\begin{align*}
\|
u\|_s^2 
& = \|u_1\|_s^2 +\|u_2\|_s^2 \\
& \le c( \|g_1\|_s^2 +\|g_2\|_s^2 )\\
& \le c(\|f_1\|_s^2 +\|f_2\|_s^2)\\
& =c\|f\|_s^2 . \qquad \qquad \qquad \qquad \qed
\end{align*}
\renewcommand{\qed}{}\end{pf}

We let $N_s$ be the operator on $\wsq$ defined by
\begin{equation}\label{box}
(\dbar\dbar^* +\dbar^*\dbar)N_s f=f   
\end{equation}
for all $f\in\wsq$.  [Notice that the harmonic space for 
the operator on the left side of (\ref{box}) is just the 
zero space---by the preceding arguments.  Therefore
this last condition uniquely defines $N_s$.]
We call $N_s$ the {\em Neumann operator} for the
$\dbars$-Neumann problem.
Thus we have proved that $N_s$ is a bounded operator from
$\wsq$ into itself, for $0<q\le n$.

\begin{remark}{\em 
We want to stress the fact that the results of the present section 
are {\em independent} of the particular choice of
the Sobolev inner product.  
In fact, the same arguments work for any $s\ge0$, not necessarily
integral, and any choice of an equivalent norm in $\wsq$. 
The form of the inner product (\ref{Sobolev}) will only become
relevant in the next section 
since we seek explicit formulas for $\dbar^*$ and its domain.
}
\end{remark}

\section{The Hilbert space adjoint $\dbar^*$ of $\dbar$}

In this section we wish to make the $\dbars$-Neumann problem more
explicit by calculating the domain of $\dbar^*$, and the operator
$\dbar^*$ itself by showing how
it actually operates on the $(0,q)$-forms in its domain.
As a result,
we shall formulate  problem (\ref{dbar-neumann})
as a boundary value problem in which
the equation on the domain is of the form $\Box+G_s$ where $\Box$
is the complex Laplacian and $G_s$ is a so-called {\em singular
Green's operator} (see \cite {Grubb}). This result 
demonstrates a striking difference with
the classical case of the $(\dbar,0)$-Neumann problem, where the
adjoint 
is taken with respect to the $L^2$-inner product, and no operator
$G_s$ appears.  
The particular expression of the Hilbert space adjoint $\dbar$, and
of the singular Green's operator $G_s$ arising in the
$\dbars$-Neumann 
problem depend on the choice of the inner product in $\wsq$.  We
note that with
our definition (\ref{Sobolev}) we have 
\begin{equation}\label{iteration}
\la f,g\ra_s = \la f,g\ra_0 
	+ \sum_{j=1}^{2n} \la D_j f,D_j g\ra_{s-1} .
\end{equation}

Recall that for a $(0,q)$ form $\phi=\sum_{|J|=q}\phi_J d\bar z^J$   
with $C^\infty$ coefficients, the operator $\dbar$
is defined as
$$
\dbar\phi =\sum_{|K|=q+1} \sum_{kJ}\e{K}{kJ} \pd{\phi_J}{\bar z_k} 
d\bar z^K .
$$
Then, the {\it formal adjoint}
$\vartheta$ of $\dbar$ is easily calculated to be
$$
\vartheta \phi
= -\sum_{|I|=q-1} \e{J}{iI} \pd{\phi_J}{z_i} d\bar z^I.
$$

We want to compute the Hilbert space adjoint $\dbar^*$, together
with its domain. In carrying out this program, a central role is
played by a particular extension of the normal vector field on
$b\Omega$ to a suitable tubular neighborhood of $b\Omega$.  The
differential condition arising in the description of $\dom\dbar^*$ is
most easily expressed if we make the following choices.

Let $\varrho$ be the signed Euclidean distance from $b\Omega$
(negative inside, positive outside).  In a
suitable tubular neighborhood $U$ of $b\Omega$, this function
$\varrho$ is well defined and in 
$C^\infty(\overline{U})$.  We define
the vector field $N$ on $U$ by setting $N=\text{grad\,}\varrho$.
Then $N$ is the outward unit vector field, and if we set
$N=\sum_{j=1}^{2n}\nu_j D_j$, then $N\nu_j =0$ on $U$. 
[This last one
is in fact the property that at several stages makes
our computations easier, and the formulas appearing simplier.]

\begin{prop}\label{domain} \sl 
Let $\Omega$ be a smoothly bounded domain
in $\cn$. Let $\dbar^*$ be the $\ws$-Hilbert space adjoint of
$\dbar$. Then 
$$
\dom\dbar^* \cap C_{(0,q+1)}^\infty(\bar\Omega)
=\{ \psi: N^s \bigl(\psi \llcorner\dbar\varrho\bigr)_I 
=0 \text{ on }b\Omega, \text{ for all }I, 
|I|=q \} .
$$
\end{prop}

\noindent Notice that $N^s$ denotes the $s$-fold composition of $N$ 
with itself.

Here the contraction of a $(0,q+1)$ form $\psi$ with a $(0,1)$ form 
$\omega=\sum_{k}\omega_k d\bar z_k$
is defined by the formula
$$
\psi\llcorner\omega:= \sum_I \sum_{kK}\e{K}{kI}\psi_K \bar\omega_k
d\bar z^I .
$$

\begin{remark}{\em 
Suppose that $r$ is a generic defining 
function for $\Omega$, $C^\infty$ in a neighborhood of
$\overline{\Omega}$, and we wish to express $\dom\dbar^*$ in
 terms of
$r$ and $\p/\p r$.  Then we obtain the following description.
There exists a differential operator $L_s$ of order $s$, with
$C^\infty (\overline{U})$ coefficients, whose leading
term is $(\p/\p r)^s$, and such that
$$
\dom\dbar^* \cap C_{(0,q+1)}^\infty(\bar\Omega)
= \biggl \{ \psi: L_s \bigl(\psi \llcorner\dbar r\bigr)_I 
=0 \text{ on }b\Omega, \text{ for all }I, |I|=q \biggr \} .
$$
Indeed, on $U$,
$$
\pd{}{r} = |\text{grad\,} r|N+ gX,
$$
where $g\in C^\infty(\overline{U})$, $g=0$ on $b\Omega$, and $X$ is
a vector field on $U$.
}
\end{remark}

We set $\dbar^* =\vartheta+\K$, and we 
want to determine $\K$. Our result is
the following.

\begin{prop}\label{K} \sl
Let $\Omega$ be a smoothly bounded domain in $\cn$, and 
let $s$ be a positive
integer. Let $\dbar^*$ 
be the $\ws$-Hilbert adjoint of $\dbar$, 
and $\vartheta$ be the formal adjoint
of $\dbar$, 
respectively. Set $\dbar^*=\vartheta + {\cal K}$. Then, for
a $(0,q+1)$ form $\psi$, $\K\psi:=\sum_{|I|=q}\omega_I d\bar z^I$ 
is the $(0,q)$ form whose components are solutions of the following 
\bvp:
\begin{equation}\label{BV}
\begin{cases}
\sum_{j=0}^{s}(-\Delta)^j \omega_I =0 & \text{ on }\Omega\\
\sum_{j=0}^{s+k-1}T_j N^{s+k-1-j}\omega_I =P_{s+k}^{(I)}\psi  
	& \text{ on }b\Omega ,\, k=1,\dots,s \, .
\end{cases} 
\end{equation}
Here $T_k$ denotes a
tangential differential operator of order $\le k$, with $C^\infty$
coefficients, $T_0 =(-1)^{k-1} \cdot \hbox{id}$ 
on $b\Omega$, and $P_{s+k}^{(I)}$ is
a differential operator of order $s+k$ with $C^\infty$ coefficients 
and acting on the components of $\psi$.
\end{prop}

As a consequence, from the
theory of elliptic boundary value problems \cite{LiMa}, 
we shall obtain the next
result. 
\begin{cor}\label{K-order-1}
Let $s$ a be positive integer.  Then, $\K$ is a
well defined operator of order $1$. More
precisely, for all $t>s+1/2$.  
there exists a positive constant $C_t >0$ 
such that we have the estimate
$$
\|\K \psi\|_{t-1} \le C_t \|\psi\|_t
$$
for all $\psi\in C^\infty_{(0,q+1)}(\overline{\Omega})$.
Furthermore, when restricted to purely tangential forms,
$\K$ is of order $0$, i.e. for all $t>s+1/2$ 
there exists $C_t >0$ such that 
if $\psi\llcorner\dbar\varrho=0$ in a neighborhood of 
$b\Omega$, then
$$
\|\K\psi\|_{t-1} \le C_t \|\psi\|_{t-1} .
$$ 
\end{cor}
\noindent As a consequence of these facts,
we obtain the following representation for the
$\dbars$-Neumann problem. 
We set $G_s :=\dbar\K+\K\dbar$.
With the notation above, the $\dbars$-
Neumann problem is equivalent
to the boundary value problem
$$
\begin{cases}
(\Box+G_s)u=f &\quad\text{on }\Omega\\
N^s (u\llcorner\dbar\varrho)=0 &\quad\text{on }b\Omega\\
N^s (\dbar u\llcorner\dbar\varrho)=0 &\quad\text{on }b\Omega \, . 
\end{cases}
$$

Here $\Box
:=\dbar\vt+\vt\dbar$ is the complex Laplacian,
and it equals $-4\Delta$ on $\Omega\ss{\Bbb C}^n$. 
Notice that $G_s$ is the singular Green's operator we mentioned
earlier. The operator $G_s$ is of order $2$, so of the same order as 
the complex Laplacian $\Box$. Moreover notice that $G_s u$
only depends on the boundary values of $u$ and $\dbar u$ and their
derivatives up to order $2s$, and that in
general $G_s$ is not diagonal. An analysis of the analogue of the
operator $G_s$ in the case of the de Rham complex, appears in
\cite{FKP1}. 
\begin{pf*}{Proof of Theorem 3.1}
Let $\phi\in C^\infty_{(0,q)} (\overline{\Omega})$ and
$\psi\in C^\infty_{(0,q+1)} (\overline{\Omega})$.
Using Green's formula we have
\begin{align}\label{>>>}
\la \dbar\phi,\psi\ra_s
= & \la \phi,\dbar^* \psi\ra_s = \la \phi,\vartheta \psi\ra_s
	+\la\phi,\K\psi\ra_s \notag \\
= & \la \phi,\vartheta \psi\ra_s + \sum_{0\le|\alpha|\le s} \CP
	\sum_{KkI} \e{K}{kI}
\int_{b\Omega}D^\alpha \phi_I \overline{D^\alpha \psi_K}
\pd{\varrho}{\bar z_k} .
\end{align}
Recall that the $(0,q+1)$ form $\psi$
belongs to $\dom \dbar^*$ if and only
if there exists a 
constant $C_\psi >0$ such that 
$|\la \dbar \phi,\psi\ra_s|\le C_\psi \|\phi\|_s$  
for all $\phi\in\dom\dbar$. 
Hence $\psi\in\dom\dbar^*$ if and only if the
boundary terms in the calculation (\ref{>>>}) above 
can be bounded by $C_\psi \|\phi\|_s$. 
By the Sobolev trace theorem we can bound the terms of the form
$$
\int_{b\Omega} D^\alpha \phi_I \overline{D^\alpha \psi
_K}
\pd{\varrho}{\bar z_k} 
$$ 
when $|\alpha|\le s-1$. Thus it suffices to consider the sum
$$
\sum_{|\alpha|= s}
\sum_{KkI}
\int_{b\Omega}D^\alpha \phi_I \overline{D^\alpha \psi_K}
\pd{\varrho}{\bar z_k} .
$$
By integrating by parts we can move tangential derivatives from
$\phi$ to $\psi$, so only the $s$ normal derivatives on $\phi$ may
cause trouble. 

We decompose the standard derivatives in the coordinate directions 
into their normal and tangential components:
$$
D_j = Y_j +\nu_j N,
$$
where $N$ is the normal derivative, and $Y_j$ are tangential vector 
fields. Then
$$
D^\alpha=(Y_{\alpha_{p_1}}+\nu_{\alpha_{p_1}}N)\cdots
(Y_{\alpha_{p_s}}+\nu_{\alpha_{p_s}}N).
$$
Notice that, since $\sum_j \nu_j^2 \equiv 1$ and 
$N=\sum_j \nu_j D_j$,  
we have that $\sum_j \nu_j Y_j =0$. 
Therefore, when considering $s$ normal derivatives on
$\phi_I$, we have 
\begin{align*}
\lefteqn{\sum_{|\alpha|=s}\CP \sum_{KkI} } \\
& \biggl [ \e{K}{kI} \int_{b\Omega} 
(\nu_{\alpha_{p_1}}N)\cdots
(\nu
_{\alpha_{p_s}}N)\phi_I
\overline{(Y_{\alpha_{p_1}}+\nu_{\alpha_{p_1}}N)\cdots 
(Y_{\alpha_{p_s}}+\nu_{\alpha_{p_s}}N) \psi_K} 
\pd{\varrho}{\bar z_k} \biggr ] \\
& = \sum_{|\alpha|=s} \CP \sum_{KkI} \e{K}{kI} \int_{b\Omega} 
(\nu_{\alpha_{p_1}})^2 \cdots (\nu_{\alpha_{p_s}})^2 
\bigl( N^s \phi_I\bigr) \overline{\bigl(N^s \psi_K)}
	\pd{\varrho}{\bar z_k} \\
& = \bigl( \sum_{|\alpha|=s}\CP (\nu_{\alpha_{p_1}})^2 \cdots 
(\nu_{\alpha_{p_s}})^2 \bigr)\sum_I \int_{b\Omega} 
(N^s \phi_I) \overline{\bigl
( \sum_{Kk} \e{K}{kI} N^s (\psi_K
\pd{\varrho}{z_k} )\bigr) } . 
\end{align*}
Now, if $\psi\in C^\infty_{(0,q+1)} (\overline{\Omega})$ and
$$
0=\sum_{Kk} \e{K}{kI} N^s (\psi_K
\pd{\varrho}{z_k}) = N^s (\psi\llcorner \dbar\varrho)_I
$$
on $b\Omega$ for all $I$, then clearly $\psi\in\dom\dbar^*$.  

On the other hand, suppose that 
$N^s (\psi\llcorner\dbar\varrho)_I \neq 0$ on
$b\Omega$ for a certain $I$.  We may assume that
$$ 
\text{Re} \bigl( N^s (\psi\llcorner\dbar\varrho)_I\bigr)  \ge 1  
\quad\text{on } B(p,\delta)\cap\overline{\Omega}, 
$$
where $B(p,\delta)$ is a small
ball center at $p\in\Omega$.  For $\varepsilon>0$,
consider the collection of $(0,q)$
forms $\phi^{(\varepsilon)}$,
$$
\phi^{(\varepsilon)}:=(-\varrho)^{s-1}
	(-\varrho+\varepsilon)^{3/4} \chi d\bar z^I ,
$$
where $\chi$ is a non-negative $C^\infty$ cut-off function,
$\text{supp}\chi\ss B(p,\delta)$, and $\chi=1$ on $B(p,\delta/2)$.
Now, an easy calculation shows that
$$
\| \phi^{(\varepsilon)} \|_s \le C_1
$$
independently of $\varepsilon$, while 
$$
\left |\int_{b\Omega} N^s \phi_I^{(\varepsilon)} 
	\cdot\overline{N^s (\psi\llcorner\dbar\varrho)_I} \right | 
\ge C_2 \varepsilon^{-1/4} ,
$$
which is unbounded, as $\varepsilon\rightarrow0$.
This finishes the proof of the proposition.
\end{pf*}
\begin
{remark}{\em 
We observe that 
$\dom\dbar^* \cap C_{(0,q+1)}^\infty(\bar\Omega)$
is dense in $W^s_{(0,q+1)}(\Omega)$.  Therefore it 
suffices to show that for
any $\varepsilon>0$ and 
$\phi\in C^\infty_{(0,q+1)}(\overline{\Omega})$ there exists  
$\psi\in C^\infty_{(0,q+1)}(\overline{\Omega})$ with 
$\|\psi\|_s <\varepsilon$  and $\phi-\psi\in\dom\dbar^*$.

Having fixed $\phi$ and $\varepsilon$, let 
$\chi\in C^\infty_0 (-1,1)$ and $\chi=1$ in a neighborhood of the 
origin. Then the form $\psi$ 
$$
\psi:= (
1/s!)(-\varrho)^s \chi(-\varrho/\varepsilon) \bigl( N^s
(\phi\llcorner \dbar\varrho)\bigr) \wedge\dbar\varrho 
$$
satisfies the required conditions.
}
\end{remark}

\begin{pf*}{Proof of Proposition 3.3}
We have set $\dbar^* =\vartheta+\K$, so that for 
$\psi\in\dom\dbar^*$ we  have 
\begin{equation}\label{dag}
\la \phi,\dbar^* \psi\ra_s = \la\phi,\vartheta\psi\ra_s
+\la\phi,\K\psi\ra_s .
\end{equation}
On the other hand by (\ref{>>>}) we see that, for
$\psi\in\dom\dbar^*$ and $\phi\in C^\infty_{(0,q)
}
(\overline{\Omega})$
we have the equality
\begin{equation*}
\la\dbar\phi,\psi\ra_s = \la\phi,\vartheta\psi\ra_s 
+\sum_{0\le|\alpha|\le s} \CP
\sum_{KkJ}\e{K}{kJ}\int_{b\Omega} D^\alpha
\phi_J \overline{D^\alpha \psi_K}\pd{\varrho}{\bar z_k} \, ; 
\end{equation*}
so it follows that
\begin{equation}\label{bnry-eq-K}
\la \phi,\K\psi\ra_s = \sum_{0\le|\alpha|\le s} \CP
\sum_{KkJ}\e{K}{kJ}\int_{b\Omega} D^\alpha 
\phi_J \overline{D^\alpha \psi_K}\pd{\varrho}{\bar z_k} .
\end{equation}
By choosing $\phi$ with compact support in $\Omega$ we find that
$\K\psi$ satisfies 
\begin{align*}
0 = & \la\phi,\K\psi\ra_s \\
= & \sum_{|J|=q}\sum_{0\le|\alpha|\le s} \CP
\int_\Omega D^\alpha \phi_J \overline{D^\alpha (\K\psi)_J} \\
= & \sum_{|J|=q}\sum_{0\le|\alpha|\le s} (-1)^{|\alpha|} \CP
	\int_\Omega \phi_J \overline{D^{2\alpha} (\K\psi)_J} .
\end{align*}
Since this holds 
for all $\phi\in C^\infty_{(0,q)} (\Omega)$ with compact support in 
$\Omega$, we see that $(\K\psi)_J$ must satisfy the equation 
$$
0=\sum_{0\le|\alpha|\le s}(-1)^{|\alpha|}\CP D^{2\alpha}(\K\psi)_J 
=\sum_{j=0}^{s} (-\Delta)^j (\K\psi)_J \quad\text{on }\Omega
$$
for all $J$, which is the equation on the interior of
$\Omega$ that appears in
(\ref{BV}). 

Now we move on to consider the boundary conditions 
that $\K\psi$ must
satisfy. For $\phi\in C^\infty_{(0,q)}(\overline{\Omega})$,
by repeatedly applying Green's theorem 
to the left hand side of equation (\ref{bnry-eq-K}), 
and recalling equation (\ref{iteration}), we have
\begin{align*}
\la \phi,\K\psi\ra_s 
& = \sum_{|J|=q} \biggl ( \la \phi_J, (\K\psi)_J \ra_0
+\sum_{j=1}^{2n}\la D_j \phi_J ,D_j (\K\psi)_J \ra_{s-1} \biggr) \\
& = \sum_{|J|=q}\biggl ( \int_\Omega \phi_J \overline{(\K\psi)_J} 
+\sum_{i=1}^{2n} \int_\Omega D_i \phi_J \overline{D_i (\K\psi)_J} \\ 
& \qquad\qquad + \sum_{1\le|\beta|\le s-1} \CPb \sum_{i=1}^{2n} 
\int_\Omega D_i D^\beta \phi_J \overline{D_i 
D^\beta (\K\psi)_J} \biggr )  \\ 
& = \sum_{|J|=q} \biggl(\int_{b\Omega} \phi_J
\overline{N(\K\psi)_J} + \sum_{1\l
e|\beta|\le s-1} \CPb
\int_{b\Omega}  D^\beta \phi_J \overline{ND^\beta(\K\psi)_J} \\
& \qquad\qquad 
-  \la \phi_J , \Delta(\K\psi)_J \ra_{s-1} + \dots \biggr ) , 
\end{align*}
where the dots stand for terms that do not contribute 
to any boundary expression.

We iterate this calculation on the last term on the right in the
above chain of equalities to obtain that
$$
\la\phi,\K\psi\ra_s  = \sum_{|J|=q} 
\sum_{i=0}^{s-1} \sum_{|\alpha|\le i}\CP \int_{b\Omega}
(D^\alpha\phi_J)\overline{ND^\alpha(
-\Delta)^{s-1-i}(\K\psi)_J} +
\dots , 
$$
where the dots have the same meaning as before.
>From this equation and (\ref{bnry-eq-K}) it follows that, 
for all $J$,
\begin{multline}\label{*}
\sum_{i=0}^{s-1} \sum_{|\alpha|\le i} \CP
\int_{b\Omega}(D^\alpha\phi_J) 
\overline{ND^\alpha(-\Delta)^{s-1-i}(\K\psi)_J} \\
= \sum_{0\le|\alpha|\le s} \CP \sum_{kK} \e{K}{kJ} \int_{b\Omega}
D^\alpha \phi_J \overline{D^\alpha \psi_K} \pd{\varrho}{\bar z_k} . 
\end{multline}
This equation must hold true for all
 $\phi\in C^\infty_{(0,q)}
(\overline{\Omega})$. Thus we need to isolate the terms containing 
$N^\ell \phi_J$ for $\ell=0,1,\dots,s-1$, and for all $J$.

Now observe that, if $f$ and $g$ are smooth functions on the
boundary, then 
\begin{align*}
\sum_{j=1}^{2n} \int_{b\Omega} D_j f \overline{D
_j g}
= & \sum_j \int_{b\Omega}(Y_j +\nu_j N)f
			\overline{(Y_j +\nu_j N)g}\\
= & \sum_j \int_{b\Omega} Y_j f \overline{Y_j g}
		+\int_{b\Omega} Nf \overline{N g}, 
\end{align*}
where we have used the fact that $\sum_j \nu_j Y_j =0$. Now 
$$
D^\alpha = T_{\alpha,|\alpha|} +T_{\alpha,|\alpha|-1}N+\cdots
+\nu^\alpha N^{|\alpha|},
$$
where $T_{\alpha,k}$ is a tangential operator of order 
$\le k$, and $\nu:=$ $(\nu_1,\dots,\nu_{2n})$. 
Therefore the left hand side of (\ref{*}) equals
\begin{align}
\lefteqn{
\sum_{i=0}^{s-1} \sum_{|\alpha|\le i} \CP \int_{b\Omega} \bigl( 
T_{\alpha,|\alpha|}+T_{\alpha,|\alpha|-1}N+\cdots+\nu^\alpha
N^{|\alpha|} \bigr)\phi_J} \\
& \qquad \qquad \cdot\overline{N D^\alpha (-\Delta)^{s-1-i}
(\K\psi)_J }\notag \\
& = \sum_{\ell=0}^{s-1} \biggl( \int_{b\Omega} 
N^\ell \phi_J \cdot \overline{ \bigl[ \sum_{i=\ell}^{s-1} 
\sum_{\ell\le|\alpha|\le i} \CP T_{\alpha,|\alpha|-\ell}^* 
ND^\alpha (-\Delta)^{s-1-i} (\K\psi)_J \bigr]} \biggr) \notag\\
& = \sum_{\ell=0}^{s-1} \biggl( \int_{b\Omega} 
N^\ell \phi_J \cdot \overline{ \sum_{\ell\le|\alpha|\le s-1} 
\bigl[ \sum_{j=0}^{s-1-|\alpha|} \CP T_{\alpha,|\alpha|-\ell}^* 
ND^\alpha (-\Delta)^{j} (\K\psi)_J \bigr]} \biggr) .
\label
{DAG} 
\end{align}
Notice that in the above calculations we have obtained the identity  
\begin{multline}\label{K-identity}
\la\phi,\K\psi\ra_s = \sum_J \biggl( \la\phi_J, 
\sum_{j=0}^{s}  (-\Delta)^j (\K\psi)_J \ra_0 \\
+ \sum_{\ell=0}^{s-1}  \int_{b\Omega} 
N^\ell \phi_J \cdot \overline{ \sum_{\ell\le|\alpha|\le s-1} 
\bigl[ \sum_{j=0}^{s-1-|\alpha|} \CP T_{\alpha,|\alpha|-\ell}^* 
ND^\alpha (-\Delta)^{j} (\K\psi)_J \bigr]} \biggr) .
\end{multline}
In particular, for $\nu:=$ $(\nu_1,\dots,\nu_{2n})$,
we have that $T_{\alpha,0}=\nu^\alpha=T_{\alpha,0}^*$ and for $\ell$ 
a positive integer we have
\begin{equation}\label{***}
\sum_{|\alpha|=\ell-1}\CP \nu^\alpha D^\alpha=
\sum_{|\beta|=\ell-2}\CPb \nu^\beta \bigl(\sum_{i=1}^{2n}\nu_i
D_i\bigr)D^\beta=\dots= N^{\ell-1}.
\end{equation}
Thus the last summand on the right hand side of (\ref{DAG})
(corresponding to $\ell=s-1$)
becomes 
\begin{equation*}
\int_{b\Omega} N^{s-1} \phi_J\cdot \overline{\biggl(
\sum_{|\alpha|=s-1} \CP
T_{\alpha,0}^* \bigl[ ND^\alpha (\K \psi)_J \bigr] \biggr)} 
= \int_{b\Omega} N^{s-1} \phi_J\cdot \overline{N^s (\K\psi)_J} .
\end{equation*}
The right hand side of (\ref{*}) can be treated in the same way:
\begin{multline}
\sum_{0\le|\alpha|\le s} \CP \int_{b\Omega} \biggl( 
\sum_{\ell=0}^{|\alpha|}
T_{\alpha,|\alpha|-\ell} N^\ell \phi_J\biggr)\overline{\biggl(
\sum_{kK} \e{K}{kJ} D^\alpha \psi_K \pd{\varrho}{z_k} \biggr)}\\  
= \sum_{\ell=0}^{s} \int_{b\Omega} N^\ell \phi_J \cdot 
\overline{ \sum_{\ell
\le|\alpha|\le s} \CP
T_{\alpha,|\alpha|-\ell}^* \biggl(
\sum
_{kK} \e{K}{kJ} D^\alpha \psi_K \pd{\varrho}{z_k} \biggr)} .
\label{DDAG}
\end{multline}
Notice that the top order term vanishes since $N^s \phi_J$ is
paired with
$$
\sum_{|\alpha|=s} \CP 
T_{\alpha,0}^* \biggl( \sum_{kK} \e{K}{kJ} D^\alpha
\psi_K \pd{\varrho}{z_k} \biggr) = \sum_{kK}\e{K}{kJ} N^s \psi_K 
\pd{\varrho}{z_k} ,
$$
which equals $0$ on $b\Omega$, because $\psi\in\dom\dbar^*$. 
>From these calculations, and by equating the
 right hand sides of
(\ref{DAG}) and (\ref{DDAG}),
we obtain the $s$ boundary equations. Set 
$$
\sum_{kK}\e{K}{kJ}D^\alpha\psi_K \pd
{\varrho}{z_k} 
= (L_\alpha \psi)_J .
$$
Then, on $b\Omega$, we have
\begin{align*} 
N^s (\K\psi)_J 
& = \sum_{s-1\le|\alpha|\le s} \CP T_{\alpha,|\alpha|-s+1}^*
	(L_\alpha \psi)_J \\
\lefteqn{\sum_{s-2\le|\alpha|\le s-1} \CP 
T_{\alpha,|\alpha|-s+2}^* ND^\alpha
\biggl( 
\sum_{j=0}^{s-1-|\alpha|}(-\Delta)^{j} (\K\psi)_J \biggr)
}  
  \hbox{\qquad \qquad \qquad \qquad \qquad} \\ 
& \qquad = \sum_{s-2\le|\alpha|\le s} \CP
T_{\alpha,|\alpha|-s+2}^* (L_\alpha \psi)_J
\end{align*}
$$
\qquad \qquad \qquad \qquad \cdots \qquad \qquad  \qquad \qquad
\cdots \\  
$$
\begin{align*}
\lefteqn{\sum_{0\le|\alpha|\le s-1} \CP T_{\alpha,|\alpha|}^*
ND^\alpha 
\biggl( 
\sum_{j=0}^{s-1-|\alpha|}(-\Delta)^{j} (\K\psi)_J \biggr)} 
\hbox{\qquad \qquad \qquad \qquad \qquad} \\
& = \sum_{0\le|\alpha|\le s} \CP 
T_{\alpha,|\alpha|}^* (L_\alpha
 \psi)_J .
\end{align*}
Thus we have $s$ boundary equations in $(\K\psi)_J$. Notice that
the $k^{\rm th}$ equation has order $s+k-1$ in the normal direction,
for $k=1,\dots,s$. Since $T_{\alpha,0}^* =\nu^\alpha$ and
$-\Delta=-N^2 +T_1 N+T_2$, using formula (\ref{***}), 
the operator on
the left hand side in the $k^{\rm th}$ 
equation becomes 
\begin{align*}
\lefteqn{
\sum_{s-k\le|\alpha|\le s-1} \CP 
T_{\alpha,|\alpha|-s+k}^* ND^\alpha
\biggl( \sum_{j=0}^{s-1-|\alpha|}(-\Delta)^{j}  \biggr)}\\ 
& = N^{s-k+1} \sum_{j=0}^{k-1}(-\Delta)^{j} + \cdots +
\sum_{|\alpha|=s-1} \CP T_{\alpha,k-1}ND^\alpha \\
& = (-1)^{(k-1)}N^{s+k-1} +T_1 N^{s+k-2}+\cdots+ T_{s+k-2}N 
\end{align*}  
as in the statement of the proposition,
while the right hand side in the same equation is an operator of 
order $s+k$ (one order larger than the left hand side), that we
denote by $P^{(J)}_{s+k}$.  Then we have
\begin{equation}\label{P-s+k}
P^{(J)}_{s+k} (\psi) = 
\sum_{s-k\le|\alpha|\le s} \CP 
T_{\alpha,|\alpha|-s+k}^* (L_\alpha
 \psi)_J .
\end{equation}
This finishes the proof.
\end{pf*}

Before proving Corollary \ref{K-order-1} we need one more result.
Consider the \bvp\ (\ref{BV}) that defines the components of $\K$:   
\begin{equation}\label{BVP2}
\begin{cases}
\sum_{j=0}^{s}(-\Delta)^j u =0 & \text{ on }\Omega\\
\sum_{j=0}^{s+\ell}T_j N^{s+\ell-j}u = g_\ell
	& \text{ on }b\Omega ,\, \ell=0,\dots,s-1 \, .
\end{cases} 
\end{equation}
for given $g_\ell \in C^\infty (b\Omega)$, $\ell=0,\dots,s-1$.  
Notice that the operator $\K$ applied to a form $\psi$
gives rise to the 
composition of a (non-diagonal) differential operator acting on the
components of $\psi$, the restriction to the boundary $b\Omega$, and 
the solution operator $S$ of the (scalar) \bvp\ (\ref{BVP2}).
Then we have the following.
\begin{lem}\label{ellipticity}
The \bvp\ (\ref{BVP2}) is an elliptic 
\bvp\ with trivial kernel, that
is if $g_\ell =0$ for $\ell=0,\dots,s-1$, then 
$S(g_0,\dots,g_{s-1})=0$. 
\end{lem}
\begin{pf}
In order to prove that the \bvp\ (\ref{BVP2})
is elliptic, we use the
standard definition, see (10.1.1) in \cite{HO2}.
Given any point $p\in b\Omega$ we need to consider a $C^\infty$
change of coordinates that takes $p$ into the origin, flattens the  
boundary, and such that the transformed vector fields at the 
origin coincide with the new basis vector fields.  
We write the new coordinates as 
$(x_0,x)\in [0,+\infty)\times{\Bbb R}^{2n-1}$. Then, the normal
vector field is $\po$, and $\p_1,\dots,\p_{2n-1}$ are the tangential 
vector fileds.
After taking the Fourier transform in the tangential directions,
writing $\xi\in {\Bbb R}^{2n-1}$ for the variable dual to $x$, 
we need to show that the ordinary differential equation
\begin{equation}\label{ODE}
\begin{cases}
(-\po^2 +|\xi|^2)^s v & =0 \quad \text{on } [0,+\infty) \\
B_{s,\ell}\, v (0) & =0 \quad \ell=0,1,\dots,s-1
\end{cases}
\end{equation}
admits the trivial solution as the only bounded solution on
$[0,+\infty)$.  Here $B_{s,\ell}$ denote the top order terms of the 
boundary  operators in (\ref{BVP2}) in our special chart, after
freezing the coefficients and taking the Fourier transform.

We begin by describing the differential operators that give the
initial conditions in (\ref{ODE}).  We then prove that the only
bounded solution of (\ref{ODE}) is in fact the trivial solution.

The boundary equations in (\ref{BVP2}) arise from
the 
identity (\ref{K-identity}).  
By considering forms of the type $\phi_J
d\bar z^J$ we may reduce to the case of functions.  We set
$u=(\K\psi)_J$. Consider the top
order terms in (\ref{K-identity}),  change coordinates, and
freeze the coefficients. Write $\alpha=(k,\alpha')$ and notice that 
$\CP=\binom{s-1}{k}\CPp$.  
Then $\p^\alpha =\po^k \p^{\alpha'}$. 
Notice that the top order
term in
$T_{\alpha,|\alpha|-\ell}$ equals $\p^{\alpha'}$, and that 
$T^*_{\alpha,|\alpha|-\ell}=(-1)^{|\alpha'|}\p^{\alpha'}$. Then  
we have that
\begin{align*}
B_{s,\ell}
& = \sum_{|\alpha'|=0}^{s-1-\ell} \bn{|\alpha'|+\ell}{\ell} \CPp 
(-1)^{|\alpha'|} \p^{2\alpha'} \po^{\ell+1} 
(-\Delta)^{s-1-\ell-|\alpha'|} \\
& = \sum_{j=0}^{s-1-\ell} \bn{j+\ell}{\ell}(-\Delta')^j
(-\Delta)^{s-1-\ell-j} \po^{\ell+1} ,
\end{align*}
where $\Delta'$ is the tangential Laplacian.
Now write $\Delta=\po^2 +\Delta'$.  We claim that the following
identity holds true
\begin{equation}\label{combinatorics}
\sum_{j=0}^{s-1-k-\ell} \bn{\ell+j}{j}\bn{s-j-\ell-1}{k}
= \bn{s}{\ell+k+1}.
\end{equation}
Assume the claim for now.
Then, it turns out that
$$
B_{s,\ell} = \sum_{k=0}^{s-1-k} (-1)^k \bn{s}{\ell+k+1}
|\xi|^{2(s-1-\ell-k)}\po^{\ell+2k+1} .
$$

Next, let $v=v_\xi$ be a bounded solution of (\ref{ODE}) for 
$\xi\neq0$.  Notice that
$v=\bigl(\sum_{\ell=0}^{s-1}c_\ell x_0^\ell\bigr)e^{-|\xi|x_0}$. 
Let $f\in C^\infty_0 (\overline{{\Bbb R}^{2n-1}_+})$.  Then 
for any $\xi\neq0$, by
assumption  and by integrating by parts we have
\begin{align*}
0 & = -\sum_{\ell=0}^{s-1} \po^\ell f (0,\xi) \overline{B_{s,\ell} 
v_\xi (0)} +\int_0^\infty f(x_0,\xi) 
\overline{(-\po^2 +|\xi|^2)^s v_\xi} dx_0 \\
& = -\sum_{\ell=0}^{s-1} \po^\ell f (0,\xi) \overline{B_{s,\ell}
v_\xi (0)} +\sum_{j=0}^{s} (-1)^j \bn{s}{j} |\xi|^{2(s-j)}
\int_0^\infty f(x_0,\xi) \overline{\po^{2j} v_\xi } dx_0 \\
& = -
\sum_{\ell=0}^{s-1} \po^\ell f (0,\xi) \overline{B_{s,\ell}
v_\xi (0)} + |\xi|^{2s} \int_0^\infty f(x_0,\xi) 
	\overline{v_\xi} dx_0 \\ 
& \qquad -\sum_{j=1}^{s} (-1)^j \bn{s}{j} |\xi|^{2(s-j)} \bigl(
|\xi|^{2(s-j)} f(0,\xi)\po^{2j-1} v_x (0) 
+ \int_0^\infty \po f(x_0,\xi) \overline{\po^{2j-1} v_\xi } dx_0 
\bigr) \\
& =  -\sum_{\ell=1}^{s-1} \po^\ell f (0,\xi) \overline{B_{s,\ell} 
v_\xi (0)} 
+ |\xi|^{2s} \int_0^\infty f(x_0,\xi) \overline{v_\xi} dx_0 \\
& \qquad +\sum_{k=0}^{s-1} (-1)^k \bn{s}{k+
1} |\xi|^{2(s-k-1)}
\int_0^\infty \po f(x_0,\xi) \overline{\po^{2k+1} v_\xi } dx_0 .
\end{align*}
By applying integration by parts $(s-1)$ more times to the 
last term in the right
hand side above, we obtain that 
\begin{equation}\label{v-xi}
0= \sum_{j=0}^{s} \bn{s}{j} |\xi|^{2(s-j)}\int_0^\infty \po^j
f(x_0,\xi) \overline{\po^j v_\xi} dx_0  
\end{equation}
for all $\xi\neq0$.  

Now, for each $\xi\neq0$ we can pick $f$ so that
$f(\cdot,\xi) =v_\xi$.  Substituting in
(\ref{v-xi}) we obtain that
$$
\sum_{j
=0}^{s} \bn{s}{j} |\xi|^{2(s-j)} \int_0^\infty |\po^j v_\xi
(x_0)|^2 dx_0 =0,
$$ 
that is, $v_\xi=0$.  

Thus, we only need to prove the claim.  
If, for $p\ge m$ we set 
$F_k (p,m):=\sum_{j=0}^{m} \binom{k+j}{j}\binom{p-j}{m-j}$, we wish
to show that
\begin{equation}\label{claim}
F_k (p,m) =\bn{p+k+1}{m} .
\end{equation}
Observe that (\ref{claim}) holds true for $m=0,1$ and $p\ge1$, and 
for $p=m$, by direct 
computation and well known properties of binomial coefficients.
Assume the statement true for $p-1$ and all $m\le p-1$.
Since 
$$
F_k (p,m)=F_k (p-1,m)+F_k (p-1,m-1),
$$
equality (\ref{claim}) follows by induction and the equality in the 
case $m=p$.  This finishes the proof of the 
ellipticity of (\ref{ellipticity}).

Finally, if all the boundary data $g_{\ell}$ in problem
(\ref{ellipticity}) are 
identically $0$, then the only solution of the boundary value 
problem is the  trivial one. In fact, if $u$ is such a solution,
the identity (\ref{K-identity}) with $u$ in place of $\K \psi $ 
implies that $u$ is  orthogonal in the
$W^s$ sense to all $\phi \in C^{\infty} (\Omega)$, hence $u=0$.
\end{pf}

Finally, we have:
\begin{pf*}{Proof of Corollary 3.4}
Clearly, $\K$ is well defined
as composition of differential operators, restriction to the
boundary, and the operator $S$ solution of the \bvp\ in the previous
Lemma. 

Next, we use standard estimates for elliptic \bvp s, as in
\cite{LiMa} Theorem 5.1, and Lemma \ref{ellipticity}. Recall that
$P^{(I)}_{s,k}$ is a differential operator of order $s+k$, containing 
$s$ at most 
derivatives in the normal direction.  Then we see that for all
$t>s+1/2$ 
\begin{align*} 
\| \K \psi\|_{t-1} 
& \le C_t \sum_I \| (\K \psi)_I\|_{t-1}  \\
& \le C_t \sum_{I}\sum_{k=1}^{s}
	 \| P^{(I)}_{s,k} \psi\|_{W^{t-1-
(s+k-1)-1/2}(b\Omega)}\\
& \le C_t  \sum_{I}\sum_{k=1}^{s}
\| N^k \psi\|_{W^{t-k-1/2}(b\Omega)}\\
& \le C_t \sum_{I}\sum_{k=1}^{s} \|N^k \psi\|_{t-k}\\
& \le C_t \|\psi\|_t ,
\end{align*} 
where we use the assumption $t>s+1/2$ in order to able be to apply 
the trace theorem.

Finally notice that $\psi\llcorner\dbar\varrho=0$ 
in a neighborhood of $b\Omega$,
$P^{(I)}_{s,k}$ becomes an operator of one degree lower, i.e., of
order $s+k-1$.  Repeating the argument above, we obtain that, for
$t>s+1/2$
$$
\|\K\psi\|_{t-1} \le C_t \|\psi\|_{t-1} 
.
$$ 
This concludes the proof of the corollary.
\end{pf*}

{\bf Final Remarks.}
The results of Section 3 are obtained under 
a specific formulation of
the Sobolev inner product.  If we modify the
formulation by choosing
other positive coefficients 
$\gamma_\alpha$ in the definition of the inner product
(\ref{Sobolev}), results analogous to those
presented here should still hold. 
It is also the case that the formulas that arise in these 
formulations of the norm are probably much less tractable.
  
The situation seems quite different if we 
take a generic equivalent norm.
Consider, for instance, the weighted theory
of the $\dbar$-Neumann problem, as developed by Kohn in \cite{Kohn}.
Kohn showed that the regularity properties enjoyed by the canonical
solution in the weighted case are in general much stronger than the
ones enjoyed by the classical canonical solution (see also
the aforementioned work of Christ [Ch]).  Therefore, it is
clear that much has still to be understood in the general case. 
We shall
provide no details about the treatment of equivalent Sobolev
topologies.

In the present paper we have worked with $(0,q)$
forms on a domain $\Omega$ in $\cn$. These results hold true in
the case of $(p,q)$ forms, with no change in the proofs. 
Routine modifications (see \cite{Folland-Kohn}) should allow one 
to work out the case of a smoothly bounded pseudoconvex domain $M'$
in a complex, or even an almost complex, manifold $M$. 

Of course it is also of interest to work out sharp estimates for
the $\dbars$ problem, and to calculate the full Hodge and
spectral theories; we save that work for a future series of
papers.

\end{document}